\documentclass[11pt]{article}

\usepackage{graphicx}
\usepackage{amsmath,amssymb}
\usepackage{mathrsfs}
\usepackage{booktabs,multirow}

\newtheorem{theorem}{Theorem}
\newtheorem{claim}{Claim}
\newtheorem{lemma}{Lemma}

\newtheorem{definition}{Definition}

\def\square{\hbox{\vrule height8pt depth0pt
\vbox{\hrule width7.2pt\vskip7.2pt\hrule width7.2pt}\vrule
height8pt depth0pt}\smallskip}

\title{\bf On some path-critical Ramsey numbers}

\author{
Ye Wang,\footnote{College of Mathematical Sciences, Harbin Engineering University, Harbin 150001, China. Email: {\tt ywang@hrbeu.edu.cn}. Supported in part by NSF of Heilongjiang Province of China (No.\ LH2021A004).}\;\;\;\;
Yanyan Song\footnote{Corresponding author. College of Mathematical Sciences, Harbin Engineering University, Harbin, 150001, China. Email: {\tt songyanyan@hrbeu.edu.cn}.}
}

\begin{document}
\date{}
\maketitle
\begin{abstract}
For graphs $G$ and $H$, the Ramsey number $R(G,H)$ is the smallest $r$ such that any red-blue edge coloring of $K_r$ contains a red $G$ or a blue $H$.
The path-critical Ramsey number	$R_{\pi}(G,H)$ is the largest $n$ such that any red-blue edge coloring of $K_r \setminus P_{n}$ contains a red $G$ or a blue $H$, where $r=R(G,H)$ and $P_{n}$ is a path of order $n$.
In this note, we show a general upper bound for $R_{\pi}(G,H)$, and determine the exact values for some cases of $R_{\pi}(G,H)$.

\medskip
	
\end{abstract}

\section{Introduction}
\indent
\par For graphs $F$, $G$ and $H$, let $F\rightarrow (G,H)$ signify that any red-blue edge coloring of $F$ contains a red $G$ or a blue $H$. The Ramsey number  $R(G,H)$ is defined as the smallest integer $r$ such that $K_{r}\rightarrow (G,H)$. For a graph $F$ and a positive integer $r$, if $|V(F)|\le r$, we denote by $K_r\setminus F$ the graph obtained from $K_r$ by deleting the edges of $F$ from $K_r$.

One significant issue is to determine the ``largest" possible graph $F$ such that $K_{r}\setminus F \rightarrow (G,H)$. If we consider ``largest" in terms of the largest number of edges, this problem is related to the size Ramsey number, which was defined by Erd\H{o}s, Faudree, Rousseau and Schelp \cite{Erdos 1978} as
$$\hat{r}(G, H)=\min \{|E(F)|: F \rightarrow (G, H)\}.$$
If the ``largest" possible graph $F$ is a star, then this problem is associated with the star-critical Ramsey number defined by Hook and Isaak \cite{Hook 2011}.
Wang and Li \cite{Wang 2020} extended the definition to restrict $F$ to the matching, path and complete graph, and defined the matching-critical Ramsey number, path-critical Ramsey number and complete-critical Ramsey number. For results on star-critical Ramsey numbers and related critical Ramsey numbers, see a survey \cite{Budden 2023}.

In this note, we consider the path-critical Ramsey number, which is defined as follows, where $P_n$ is a path on $n$ vertices.
\begin{definition}
For graphs $G$ and $H$, define the path-critical Ramsey number as
$$
R_{\pi}(G,H)=\max \left\{n:K_r \setminus P_n \rightarrow(G, H) \right\},
$$
where $r=R(G,H)$.
\end{definition}
Similarly, we define $R_{\mu}(G,H)$ the matching-critical Ramsey number and $R_{\omega}(G,H)$ the complete-critical Ramsey number.

Denoted by $\delta(G)$, $\Delta(G)$ and $ \chi(G)$ the minimum degree, the maximum degree and the chromatic number of $G$, respectively.
The minimum number of vertices in a color class among all proper vertex colorings of $G$ by $ \chi(G)$ colors is denoted by $s(G)$, and call $s(G)$ the chromatic surplus.
For connected graph $G$ and graph $H$ with $|V(G)|\ge s(H)$,
 Burr \cite{Burr 1981} showed that
\begin{equation*}
R(G,H)\ge(\chi(H)-1)(|V(G)|-1)+s(H).
\end{equation*}
 If the equality holds, $G$ is called $H$-good.
Let $G+H$ be the disjoint union of $G$ and $H$ by adding complete edges between $V(G)$ and $V(H)$, $G \cup H$ the disjoint union of $G$ and $H$, and $nG$ the $n$ vertex-disjoint copies of $G$.

In this note, we shall give a general upper bound for $R_{\pi}(G,H)$ and show the exact values of some cases of $R_{\pi}(G,H)$. In particular, the general upper bound is sharp for cases in $R_{\pi}(K_{1,n},K_m)$, $R_{\pi}(K_{1,n},B_m)$ and $R_{\pi}(F_{n},K_3)$, where $B_m=K_2+mK_1$ and $F_n=K_1+nK_2$.

\begin{theorem}\label{l1}
For connected graph $G$ of order $n$ and graph $H$ with $\Delta(G)=n-1$ and $n \geq s(H)$,
if $r=R(G, H) \leq(\chi(H)-1)n+s(H)-1$, then
$$
R_{\pi}(G,H)\leq \big[r-(\chi(H)-1)(n-1)-s(H)+1\big]n-1.
$$
In particular, if $G$ is $H$-good,
then $R_{\pi}(G,H)\leq n-1.$
\end{theorem}

It is shown in \cite{Wang 2020} that $K_r\setminus K_n\to (K_{1,n},K_m)$, where $r=R(K_{1,n},K_m)=n(m-1)+1$. Along with Theorem \ref{l1}, we obtain that for $m, n\geq2$,

\[
R_{\pi}(K_{1,n},K_m)=n.
\]
It is also easy to see that
$$
R_{\pi}(K_{1,m}, K_{1,n})=\left\{\begin{array}{cl}0 & \text { if $m$ and $n$ are both even, } \\ m+n-1 & \text { otherwise. }\end{array}\right.
$$

\begin{theorem}\label{t2}
For integers $m$ and $n$ with $n\geq3m+2$ and $m\geq2$,
then $R_{\pi}(K_{1,n},B_m)=n.$
\end{theorem}

\begin{theorem}\label{t3}
	For any integers $n\geq2$,
	then $R_{\pi}(F_{n},K_3)=2n.$
\end{theorem}
\begin{theorem} \label{t4}
For integers $m$ and $n$ with $n\geq2m+3$,
	then $R_{\pi}(K_{1,m},P_n)=n.$	
\end{theorem}
\begin{theorem}\label{t5}
For integers $m$ and $n$ with $n\geq m\geq1$ and $n\geq2$,
then $R_{\pi}(mK_2,nK_2)=2n+m-1.$
\end{theorem}

\medskip

\section{Proofs}
\indent

For $S\subseteq V(G)$, let $G[S]$ denote the subgraph induced by $S$ in $G$. Denote by $N_{H}(v)$ the neighborhood of a vertex $v \in V(G)$ that is contained in $V(H)$, and $d_{H}(v) = |N_{H}(v)|$. Similarly, $N^R_G(v)$ and $N^B_G(v)$ are the red neighborhood and blue neighborhood of $v \in V(G)$, respectively, and $d_G^R(v)=|N_G^R(v)|$, $d_G^B=|N_G^B(v)|$. In general, the vertex $v$ belongs to $V(G)$, but it may not belong to $V(G)$ in this note.

\medskip

\noindent{\bf Proof of Theorem \ref{l1}.} For simplicity of notation, let $k=\chi(H), s=s(H)$ and $t=r-(k-1)(n-1)-s+1$. Consider the graph $F=K_r \setminus P_{tn}$. For any red-blue edge coloring of $F$, the subgraphs of $F$ induced by red edges and blue edges are denoted by $R$ and $B$, respectively. Coloring the edges of $F$ such that
$$R=(k-t-1)K_{n-1}\cup t(K_{n}\setminus P_{n})\cup K_{s-1}.$$
Since $G$ is connected of order $n \geq s$ and $\Delta(K_{n}\setminus P_{n})=n-2$, none of $K_{n-1}$, $K_{n}\setminus P_{n}$ and $K_{s-1}$ contains $G$.
And $B$ contains no $H$ as $B$ is a $k$-partite graph with the chromatic surplus being $s-1$.
Thus $R_{\pi}(G,H)\leq tn-1=\big(r-(k-1)(n-1)-s+1\big)n-1. $
If $G$ is $H$-good,
then $r=(k-1)(n-1)+s$, which implies
$R_{\pi}(G,H)\leq n-1.$
\hfill $\square$
\medskip

\begin{lemma}\label{l4}(\cite{Erdos 1988})
For integers $m$ and $n$ with $m \geq 2$ and $n \geq 3 m-4$, then
$R\left(K_{1,n}, B_m\right)=2n+1$.
\end{lemma}

\noindent{\bf Proof of Theorem \ref{t2}.}~By Lemma \ref{l4} and Theorem \ref{l1}, we have $R_{\pi}(K_{1,n},B_m)\leq n.$ For the lower bound, consider the graph $G=G_1+G_2$ with
$G_1=K_{n+1}$ and $G_2=K_{n}\setminus P_{n}$. For any red-blue edge coloring of $G$, we shall show that $G$ contains neither a red $K_{1,n}$ nor a blue $B_m$.

Assume that $G$ contains a blue $K_{3}$. Suppose not, then for any vertex $v\in V(G_1)$, $d(v)= 2n$. As $G$ contains no red $K_{1,n}$, we have $d^B_{G}(v) \geq n+1$. Then $N^B_{G}(v)$ induces a red subgraph of order at least $n+1$, which must contain a red $K_{1,n}$, a contradiction.
Let the vertices of the blue $K_3$ be $u_1$, $u_2$ and $u_3$.
Then $d^B_{G}(u_i) \geq n-1$ for $i=1,2,3$.
Let $V_{1}=N^B_{G}(u_1) \setminus\{u_2, u_3\}, V_{2}=N^B_{G}(u_2) \setminus(N^B_{G}(u_1) \cup\{u_1\}), V_{3}=N^B_{G}(u_3) \setminus(N^B_{G}(u_1) \cup N^B_{G}(u_2))$. Since $G$ contains no blue $B_{m}$, we have
$$|V_{1}| \geq n-3,|V_{2}| \geq (n-3)-(m-2)=n-m-1,|V_{3}| \geq(n-3)-2(m-2)=n-2m+1.$$
Then
$$
\begin{aligned}
2n+1=|V(G)| & \geq |\{u_1, u_2, u_3\} \cup V_{1} \cup V_{2} \cup V_{3} | \\
&\geq 3+(n-3)+(n-m-1)+(n-2m+1)=3n-3m,
\end{aligned}
$$
which implies  $n \leq 3m+1$, a contradiction.\hfill $\square$
\medskip

\begin{lemma}\label{l5}(\cite{Li 1996})
For any integer $n \geq 2, R(F_n, K_3)=4 n+1$.
\end{lemma}

\begin{lemma}\label{l6}(\cite{Lorimer 1984})
For any integer $n \geq 2, R(n K_2, K_3)=2 n+1$.
\end{lemma}

The following claim is an easy observation, which can also be found in \cite{Li 2015}.

\begin{claim}\label{l7}
Let $n \geq 2$ be an integer and $r=R(n K_2, K_3)=2 n+1$. Define a class of red-blue edge coloring of $K_{2 n}$ as $\mathbb{H}=\left\{H_1, H_2, \ldots, H_{\lceil n / 2 \rceil-1}\right\}$, such that
$$
\begin{aligned}
H_i: H_i^R & =K_{2 i+1} \cup K_{2 n-2 i-1} \\
H_i^B & =K_{2 i+1,2 n-2 i-1},
\end{aligned}
$$
where $1 \leq i \leq\lceil n / 2\rceil-1$. For any $\left(n K_2, K_3\right)$-free coloring of $K_{r-1}$, the resulting graph must belong to the class $\mathbb{H}$.
\end{claim}

\begin{claim}\label{c1}
$R_{\pi}(F_{2},K_3)=4.$
\end{claim}

\noindent{\bf Proof.}
The upper bound follows from Theorem \ref{l1}.
For the lower bound, we shall prove that $G=K_{9}\setminus P_{4} \to (F_{2},K_3)$. Let $G=G_1+G_2$ with
$G_1=K_{5}$ and $G_2=P_{4},$ where $P_{4}:v_1v_2v_3v_4$. Suppose $G \not\to (F_{2},K_3)$.
Note that $d^B_{G}(x) \leq 4$ for any vertex $x \in V(G)$, otherwise we obtain a red $F_{2}$.
Let $H$ be the complete graph induced by $V(G_1)\cup \{v_2,v_3\}$.
Then for any vertex $x \in V(H)$, by Lemma \ref{l6}, we obtain $d^R_{H}(x) \leq 4$, thus $d^B_{H}(x) \geq 2$. So for any vertex $x \in V(H)$, we have $2 \leq d^B_{H}(x) \leq 4$.
If $d_H^B(x)=3$ for any vertex $x \in V(H)$, this contradicts the fact that $e(H^B)=\frac{\sum_{x \in V(H)} d_H^B(x)}{2}$.
Then we consider the following two cases.
\medskip

\noindent{\bf Case 1.} If there exists a vertex $x \in V(H)$ such that $d^B_{H}(x)=4$.
Let $H_{1}$ be a red $K_4$ induced by $N_H^B(x)$, and $H_{2}$ the graph induced by $V(G) \setminus(V(H_{1}) \cup \{x\})$. Note that any edge between $x$ and $V(H_{2})$ must be red.
For any vertex $y \in V(H_{2})$, we have $d^R_{H_{1}}(y) \leq 1$, otherwise yielding a red $F_{2}$.
So for any vertex $y_{1},y_{2} \in V(H_{2})$ with edge $y_{1}y_{2} \in E(G)$, we have $|N_{H_{1}}^B(y_1) \cap N_{H_{1}}^B(y_2)| \neq \emptyset$, then $H_{2}$ is a red graph,
yielding a red $F_{2}$ with $V(F_{2})= \{x\} \cup V(H_{2}) $, a contradiction.
\medskip

\noindent{\bf Case 2.} If there exists a vertex $x \in V(H)$ such that $d^B_{H}(x)=2$, let $ N_H^B(x)=\{u_1, u_2\}$.
Then $u_1 u_2$ is red and $G[N^R_{H}(x)]=K_{4}$ is $(2K_{2}, K_{3})$-free.
By Claim \ref{l7}, $K_{4}$ is the graph $H_{0}$ in $H$ with red clique $X=\{x_1\}$ and $Y=K_{3}$.
For any vertex $y \in \{u_1, u_2, v_1, v_4 \}$, we have $d^R_{Y}(y) \leq 1$, thus $d_{Y}^{B}(y) \geq 1$.
And if $x_1 y \in E(G)$, then $x_1 y$ is red.

\par If $x \in V(G_{1})$, let $x v_1$ and $x v_4$ be both red or blue, and then we obtain a red $F_{2}$, a contradiction. Without loss of generality, we may assume that $x v_1$ is red and $x v_4$ is blue, then $x_1=v_3$. Note that $d_{Y}^{B}(v_1) \geq 2$ and $d_{Y}^{B}(u_i) \geq 2$ for $i=1,2$.
So $u_{i}v_{1}$ is red, and $\{u_1, u_2, v_1, v_3, v_4 \}$ induces a red graph, yielding a red $F_{2}$, a contradiction.
If $x \in V(G_{2})$, without loss of generality, let $x=v_2$, and then $d_{Y}^{B}(v_4) \geq 2$.
Since $d_{Y}^{B}(u_i) \geq 2$ for $i=1,2$, $u_{i}v_{4}$ must be red, and $\{u_1, u_2, v_4, x_1\}$ induces a red clique. So $v_{1}v_{2}$ is red, otherwise we can get a red $F_{2}$.
Similarly, we have $d^R_{Y}(v_{1})=0$ and $x_{1}=v_3$, then $\{u_1, u_2, v_1, v_3, v_4 \}$ induces a red graph, a contradiction.
\hfill $\square$
\medskip

\begin{claim}\label{c2}
$R_{\pi}(F_{3},K_3)=6$.
\end{claim}

\noindent{\bf Proof.} The upper bound is trivial by Theorem \ref{l1}. For the lower bound, consider $G=G_1+G_2=K_{7}+K_{6}\setminus P_{6} \to (F_{3},K_3)$, where $P_{6}:v_1v_2\cdots v_6$. Let $H=G[V(G_1) \cup \{v_1, v_3, v_5 \}]=K_{10}$.
Analysis similar to that in the proof of Claim \ref{c1}, suppose to the contrary that $G$ contains neither a red $F_{3}$ with the central vertex in $V(G_{1})$ nor a blue $K_3$.
We see that $d^B_{G}(x) \leq 6$ for any vertex $x \in V(G)$.
For any vertex $x \in V(G_{1})$, we have $d^R_{H}(x) \leq 6$ by Lemma \ref{l6}, that is, $d^B_{H}(x) \geq 3$, and we conclude that $3 \leq d^B_{H}(x) \leq 6$ for any vertex $x \in V(G_{1})$. Let $F'=G[N^{B}_{H}(x)]$. We shall discuss each case as follows.
\medskip

\noindent{\bf Case 1.} If there exists a vertex $x \in V(G_{1})$ such that $d^B_{H}(x)=6$, then $3 \leq \mid V(F') \cap V(G_1) \mid \leq 6$.
Let $G'=G[N^{R}_{G}(x)]$, then $|V(G')|=6$.
For any vertex $ y \in V(G')$, we have $d_{F' \cap G_1}^R(y) \leq 1$, otherwise there is a red $F_{3}$ with the central vertex in $V(G_{1})$.
Then $d_{F' \cap G_1}^B(y) \geq |V(F') \cap V(G_1)|-1$,
which implies that $|N_{F' \cap G_1}^B(y_1) \cap N_{F' \cap G_1}^B(y_2)| \neq \emptyset$ for any vertex $ y_1, y_2 \in V(G')$.
So $G'$ is a red graph, and we get a red $F_{3}$ with the central vertex $x$, a contradiction.
\medskip

\noindent{\bf Case 2.} If there exists a vertex $x \in V(G_{1})$ such that $d^B_{H}(x)=5$, then $2\leq \mid \mathrm{V}(F') \cap V(G_1) \mid \leq 5$.

If $d^R_{G}(x)=d_G^B(x)=6$, the proof is similar to that of Case 1. If $|V(F')\cap V(G_1)|=2$, say $u_1, u_2$, then $G[N_G^R(x)]=K_6$ is $(3K_{2}, K_{3})$-free.
By Claim \ref{l7}, $K_{6}$ is the graph $H_{0}$ and $H_{1} $ in $H$ with red cliques $X$ and $Y$, where $|V(X)|\leq |V(Y)|$.
 We have $d_Y^R(u_i) \leq 1$ for $i=1,2$.
 So $d_X^R(u_i)=|X|$, and we get a red $F_{3}$ with the central vertex $u_i$, a contradiction.

If $d^B_{G}(x)=5$, then $d^R_{G}(x)=7$. As $\{v_2, v_4, v_6\}$ induces at least one red edge, we may assume $v_j v_k$ is red for $j, k \in$ $\{2,4,6\}$.
Since $G[N_{G}^R(x)]=K_{4}$ is $(2K_{2}, K_{3})$-free, by Claim \ref{l7}, $K_{4}$ is the graph $H_{0}$ in $H$ with red cliques $X=\{x_1\}$ and $Y=K_{3}$.
Then $d_{Y}^R(v_l)=0$ and $x_1 v_l \notin E(G)$ for $l \in\{2,4,6\}$ and $l \not\in\{ j, k\}$, which implies that $x_1 \in V(G_2)$ and $|V(F')$ $\cap V(G_1)| \geq 3$.
Let $V(G')=\{u_1, u_2, u_3\} \subseteq V(F') \cap V(G_1)$.
Note that $d_Y^B(u_i) \geq 2$ for $i=1,2,3$, otherwise, we have a red $F_{3}$ with the central vertex $u_{i}$.
So $u_1 x_1$, $u_2 x_1$ and $u_3 x_1$ are all red.
For any vertex $y \in N_G^R(x) \setminus\{x_1\}$, we have $d_{G'}^R(y) \leq 1$, otherwise there is a red $F_{3}$ with the central vertex $u_i$ for $i\in \{1,2,3\}$, a contradiction.
So $d_{G'}^B(y) \geq 2$, and $N_G^R(x) \setminus\{x_1\}$ induces a red graph, yielding a red $F_{3}$ with the central vertex $x$, a contradiction.
\medskip

\noindent{\bf Case 3.} If there exists a vertex $x \in V(G_{1})$ such that $d^B_{H}(x)=3$, then  $G[N^R_{H}(x)]=K_{6}$ is $(3K_{2}, K_{3})$-free. By Claim \ref{l7}, $G[N^R_{H}(x)]=K_{6}$ is the graph $H_{0}$ and $H_{1} $ in $H$ with red cliques $X$ and $Y$, where $|V(X)|\leq |V(Y)|$.

If $G[N^R_{H}(x)]=H_{0}$, then let $X=\{x_1\}$.
If $x_1 \in V(G_{1})$ with $d^B_{H}(x_{1})=5$, this proof is similar to that of Case 2.
If $x_1 \in V(G_{2})$, we have $|V(Y) \cap V(G_1)| \geq 3$.
Denote by $G'$ the graph induced by $V(G) \setminus(V(Y) \cup\{x\})$, then $|V(G')|=7$. For any vertex $y \in V(G')$, we have $d_{Y \cap G_{1}}^R(y) \leq 1$, otherwise there is a red $F_{3}$ with $V(F_{3})=\{x, y\} \cup V(Y)$ and the central vertex in $V(G_1)$.
Similarly, $V(G')$ is a red graph and $|V(G_1) \cap V(G' )| \geq 1$, which yields a red $F_{3}$ with the central vertex in $V(G_1)$, a contradiction.

Now we assume that $G[N^R_{H}(x)]=H_{1}$. If $d_G^R(x) \geq 7$, let $x v_i$ be red for $i \in\{2,4,6\}$. Then the edges between $v_i$ and $X, Y$ are completely blue, and we obtain a blue $K_{3}$, a contradiction.
If $d_G^R(x)=d_G^B(x)=6$, let $H'= G[N_G^B(x)]$.
We may suppose that there exists a vertex $ y \in V(H') \cap V(G_1)$.
Then $d_X^B(y) \geq 2$ and $d_Y^B(y) \geq 2$, and there is a blue $K_{3}$, a contradiction. So $|V(H')\cap V(G_1)|= \emptyset $, that is, $N_G^R(x) \subseteq V(G_1)$.
Let $u_1 \in V(Y) $, and we have $d_{H'}^R(u_1) \leq 3$, otherwise there is a red $F_{3}$ with the central vertex $u_{1}$.
So the edges between $X$ and $N_{H'}^B(u_1)$ are completely red, which implies that there is a red $F_{3}$ with $V(F_{3})=\{x\} \cup V(X) \cup N_{H'}^B(u_1)$ and the central vertex in $V(G_{1})$, a contradiction.
\medskip

\noindent{\bf Case 4.} For any vertex $x \in V(G_1)$, we have $d_H^B(x)=4$,
which implies that $d_H^R(x)=5$, $|V(F') \cap V(G_1)| \geq 1$ and $ 4\leq d_G^B(x)\leq 6 $.
Then $6 \leq d_G^R(x)\leq 8$.
Let $V(F')=\{u_1, x_1, x_2, x_3\}$ and $G'=G[N_H^R(x)]$, where $u_1 \in  V(G_1)$.
Then $d_{G'}^R(u_1)=2$, say $x_4$ and $x_5$. By Lemma \ref{l6}, we get a red $2 K_2$ with vertices in $N_H^R(x)$.
If there exists a vertex $x \in V(G_1)$ with $d_G^{R}(x)=8$,
then $\{v_2, v_4, v_6\}$ induces at least one red edge, which yields a red $F_{3}$ with the central vertex $x$, a contradiction.

Now we assume that $6 \leq d_G^{R}(x) \leq 7$ for any vertex $x \in V(G_1)$. If there exists a vertex $x \in V(G_1)$ such that $d_G^{R}(x)=6$, let $x v_i$ and $x v_j$ be blue for $i,j \in\{2,4,6\}$, then $v_i v_j, u_{1}v_i, u_{1} v_j$ are all red. Similarly, $x_4 x_5$ and $x_i x_j$ are blue for $i=\{4,5\}$ and $j =\{1,2,3\}$,
so $G$ contains a blue $K_{3}$, a contradiction.
If $d_G^{R}(x)=7$ for any vertex $x \in V(G_1)$, let $x v_i$ and $x v_j$ be all red and $x v_k$ blue for $i, j, k \in\{2,4,6\}$.
Then $v_i v_j$ is blue as there is a red $2 K_2$ with vertices in $N_H^R(x)$.
Without loss of generality, let $u_1 v_j$ be red and $u_1 v_i$ blue.
Similarly, $x_4 x_5$ and $v_j v_k$ are blue, and thus $v_i v_k$ is red and $v_j x_4, v_j x_5$ are blue, yielding a blue $K_{3}$, a contradiction.
 Hence, $v_j$ is adjacent to neither $x_4$ nor $x_5$.
Then we have $|N_H^B(x) \cap V(G_1)| \geq 2$.
Let $x_1 \in V(G_1)$. Since $v_k x_1$ is red, we have $v_i x_1$ is blue as there is a red $2 K_2 $ with vertices in $N_H^R(x_{1})$ and $v_j x_1$ is red.
 If $v_k x_2, v_k x_3 \notin E(G)$, we get $v_k x_4$ or $v_k x_5$ is red, yielding a red $F_{3}$ with the central vertex $u_{1}$, a contradiction.
 If $v_k x_2$ or $v_k x_3$ is red, say $v_k x_2$, then $x_{3}x_{4}$ and $x_{3}x_{5}$ are both blue, contradicting to our claim.
 \hfill $\square$
\medskip

\noindent{\bf Proof of Theorem \ref{t3}.} The upper bound follows from Theorem \ref{l1}. For the lower bound, consider the graph $G=G_1+G_2$ with $G_1=K_{2n+1}$ and $G_2=K_{2n}\setminus P_{2n}$, where $P_{2n}:v_1v_2\cdots v_{2n}$. Suppose that $G$ contains no blue $K_{3}$, and we shall show that $G$ contains a red $F_{n}$ with the central vertex in $V(G_{1})$ by induction on $n$.
The assertion holds for $n=3$ by Claim \ref{c2}, and now we assume that $n\ge 4$.

If $d_{G_1}^B(v_1)=2 n+1$, we are done. We may assume $d_{G_1}^B(v_1) \leq 2 n$. Denoted by the red edge $u_1 v_1$ with $u_{1} \in V(G_{1})$, and let $H$ be the graph induced by $\{v_1\}\cup V(G_1) \setminus\{u_1\}$. Similarly, by Lemma \ref{l6}, we have $d_{H}^R(u_1) \leq 2 n$, and we denote by $u_{2}$ the vertex adjacent to $u_{1}$ by a blue edge in $G_{1}$.
By induction, $V(G) \setminus\{u_1, u_2, v_1, v_2\}$ must induce a red $F=F_{n-1}$ with central vertex $ x \in V(G_1)$.
Let $ F'=G[N_G^R(x)]$ and $G'=G[N_G^B(x)]$. Then $G'$ is a red graph with $|V(G')| \leq 2n$.
\medskip

\noindent{\bf Case 1.} If $x u_1$ and $x u_2$ are both red, we have $x v_1$ is blue.
For any vertex $y \in V(G) \setminus (V(F) \cup \{u_1, u_2\})$, we have $x y$ is blue,
otherwise $y u_1$ and $y u_2$ are both blue, a contradiction.
Since $ | V(G) \setminus (V(F)$ $\cup \{u_1, u_2\})|=2n$, we have $|V(G')|=2n$.

If there exists a vertex $u_3 \in V(G') \cap V(G_1)$, then $u_1 u_3$ is blue, otherwise we obtain a red $F_{n}$ with $V(F_{n})=\{u_1\} \cup V(G')$ and the central vertex $u_{3}$.
So $u_2 u_3$ is red.
Similarly, $d_{G' \setminus \{u_3\}}^B(u_2)=2 n-1$, which implies that $d_{G' \setminus\{u_3\}}^R(u_1)=2 n-1$, and then we can get a red $F_{n-1}$ with the central vertex $u_{1}$.
Then we have $d_{F \setminus \{x\}}^B(u_1)=2 n-2$, otherwise there is a red $F_{n}$ with the central vertex $u_{1}$.
So $d_{F \setminus \{x\}}^R(u_3)=2 n-2$. Along with the red graph $G'$, we obtain a red $F_{n}$ with the central vertex $u_{3}$.

If $|V(G') \cap V(G_1)|=\emptyset$, then $V(F') \subseteq V(G_1)$, and $F'=K_{2 n}$ is $(nK_{2}, K_{3})$-free.
By Claim \ref{l7}, $K_{2n}$ is the graph $\mathbb{H}$ with red odd cliques $X$ and $Y$, where $|V(X)|\leq |V(Y)|$.
Let $x_1 \in V(X)$, $y_1 \in V(Y)$ and $H'= G[N_{G'}^B(y_1)]$.
If $|V(X)|=1$, we have $|V(Y)|=2n-1$.
Then $V(H') \leq 2n-1$, otherwise the edges between $x_{1}$ and $V(H')$ are completely red, and we have a red $F_{n}$ with $V(F_{n})=\{x_{1}\}\cup V(H')$ and the central vertex $x_{1}$.
Then $d_{H'}^{R}(x_1)=|V(H')|\leq 2 n-1$ and $d_{H'}^{B}(x_1)\geq 1$.
Denoted by the blue edge $x_1 u$ with $u \in V(H')$, and then $d_Y^{R}(u)=2 n-1$, which yields a red $F_{n}$ with $V(F_{n})=\{x, u\} \cup V(Y)$ and the central vertex in $V(G_{1})$.
If $|V(X)| \geq 3$, we have $|V(Y)|\leq 2n-3$.
Suppose that $d_{G'}^R(y_1) \geq 2 n-|V(Y)|+1$, and we can obtain a red $F_{n}$ with the central vertex $y_{1}$. If $d_{G'}^R(y_1) \leq 2 n-|V(Y)|$, then $|V(H')| \geq|V(Y)|$, and the edges between $V(X)$ and $V(H')$ are completely red, which yields a red $F_{n}$ with $V(F_{n})=\{x\} \cup V(X) \cup V(H')$ and the central vertex  $x_{1}$.
\medskip

\noindent{\bf Case 2.} If one of $x u_1$ and $x u_2$ is blue, without loss of generality, we assume that $x u_1$ is blue and $x u_2$ is red.
Note that there are at most two red edges between $x$  and $V(G) \setminus (V (F_{n-1}) \cup \{u_1, u_2 \})$.
If $w_{1}, w_{2} \in V(G_{2}) $ with $w_{1}w_{2} \notin E(G)$, then $x w_{1}$ and $xw_{2}$ are both red.
\par If $V(F')=V(G')=2 n$, let $x w_{1}$ be red. Then $u_2 w_{1}$ is blue.
Similarly, $u_1 w_{1}$ is red. Then we have $d_{G' \setminus \{u_1 \}}^{R}(w_{1})=0$, otherwise there is a red $F_{n}$ with the central vertex $u_{1}$.
Let $U=N^{B}_{G'\setminus\{u_1\}}(w_{1})$. So $d^{R}_{U}(u_2)=|V(U)| \geq 2 n-3 \geq 5$.
If $|(V(G') \setminus\{u_1\}) \cap V(G_1)| \geq 1$, we get a red $F_{n}$ with $V(F_{n})=\{u_2\} \cup V(G')$ and the central vertex in $V(G_{1})$.
If $|(V(G') \setminus\{u_1\}) \cap V(G_1)|=\emptyset$, then $F'=K_{2 n}$ is $(nK_{2}, K_{3})$-free.
By Claim \ref{l7}, $K_{2n}$ is the graph $\mathbb{H}$ with red odd cliques $X$ and $Y$, where $|V(X)|\leq |V(Y)|$.
Since there is a red $F_{n-1}$ with the central vertex $u_{1}$ in graph $G'$, we have $d_Y^B(u_1) \geq |V(Y)|-1$.
Hence, $d_Y^R(u_1)=|V(X)|$. By similar argument as in Case 1, we obtain a red $F_{n}$ with the central vertex in $V(G_{1})$.

\par If $V(F')=2 n+1$ and $V(G')=2 n-1$, let $x w_{1}$ and $xw_{2}$ be both red with $w_{1}w_{2} \notin E(G)$ and $w_{1}, w_{2} \in V(G_{2})$.
Then $u_2 w_{1}$ and $u_2 w_{2} $ are both blue.
Similarly, $u_1 w_{1}$ and $u_1 w_{2} $ are both red.
If there is a red $2K_{2}$ between $\{w_{1}, w_{2} \}$ and $V(G') \setminus \{u_1\}$, we have a red $F_{n}$ with $V(F_{n}) =\{w_{1}, w_{2} \} \cup V(G')$ and the central vertex $u_{1}$. Otherwise, we obtain
$|N_{G' \setminus\{u_1\}}^B(w_{1} ) \cup N_{G' \setminus\{u_1\}}^B(w_{2}) | \geq 2 n-3$. Since $G$ contains no blue $K_{3}$, the edges between $u_2$ and $N_{G' \setminus\{u_1\}}^B(w_{1} ) \cup N_{G' \setminus\{u_1\}}^B(w_{2})$ are all red, and we can get a red $F_{n-2}$ with the central vertex $u_{2}$.
Let $U'=G[N_{F'}^B(w_{1})]$.
If $d_{F'}^R(w_{1})=0$, then $U'$ is a red graph with $|V(U')| \geq 2 n-2$. So $d_{U'}^{R}(u_2)=|V(U')|$, and we obtain a red $F_{n-1}$ with the central vertex $u_{2}$.
Along with the red $F_{n-2}$ with the central vertex $u_{2}$ in $G'$,
we have a red $F_{n}$ with the central vertex $u_{2}$ for $n \geq 4$.
If $d_{F'}^R(w_{1}) \geq 1$, let $w_{1}w$ be red, where $w \in V(F')$.
By induction, there is a red $F_{n-1}$ with central vertex $x$, and there must exist a vertex $w' \in V(F')$ such that $ww'$ is red.
 Then $u_2 w'$ is blue. As $w_{1} w'$ or $w_{2} w'$ must be red, we have $u_2 w$ is blue.
  Then $u_1 w$ and $u_1 w'$ are both red, and we can get a red $F_{n}$ with $V(F_{n})=\{w, w'\} \cup V(G')$ and the central vertex  $u_{1}$.
\hfill $\square$
\medskip

\begin{lemma}\label{l9}(\cite{Parsons 1974})
For any integers $m,n$ with $n \geq 2m+1, R(K_{1,m}, P_n)=n$.
\end{lemma}

\begin{lemma}\label{l10}(\cite{Dirac 1952})
Every connected graph $G$ contains a path of length at least $\min \{2 \delta, |V(G)|-1\}$.
 \end{lemma}

\noindent{\bf Proof of Theorem \ref{t4}.} We shall prove that $G=K_{n} \setminus P_{n} \rightarrow (K_{1,m},P_n)$ for $n\ge 2m+3$ by induction on $n$.

For $n=2m+3$, we shall show that $G=K_{2m+3}\setminus P_{2m+3} \rightarrow (K_{1,m},P_{2m+3})$. Suppose that $G$ contains no red $K_{1,m}$, then
$d^{R}_{G}(v) \leq m-1$ for any vertex $v \in V(G)$, which implies $d^{B}_{G}(v) \geq m+1$. By Lemma \ref{l10}, $G$ must contain a blue $P_{2m+3}$.
For $n \ge 2m+4$, let $P_n: v_1v_2 \ldots v_n$.
By induction, $G \setminus \{v_{1} \} \rightarrow (K_{1,m},P_{n-1})$.
Then $G$ contains a blue $H=P_{n-1}$. Assume that $G$ contains no blue $P_n$, then $d^{R}_{H}(v_{1}) \geq \lceil\frac{n-1}{2}\rceil-1\geq m$, which yields a red $K_{1,m}$, completing the proof.
\hfill $\square$
\medskip

\begin{lemma}\label{l13}(\cite{Cockayne 1975, Cockayne 75,Lorimer 1984})
For any integers $m,n$ with $n \geq m\geq1, R(mK_2, nK_2)=2n+m-1$.
\end{lemma}

\noindent{\bf Proof of Theorem \ref{t5}.} By Lemma \ref{l13}, we shall prove that
$G=K_{2n+m-1}\setminus P_{2n+m-1}\to (mK_2,nK_2)$ by induction.

If $m=1$ and $n\geq 2$, it is trivial as $ K_{2n} \setminus P_{2n} \rightarrow (K_2, nK_2)$. If $m=n=2$, we shall show that $ K_{5} \setminus P_{5} \rightarrow (2K_2,2K_2)$. For the path $P_5:v_1v_2v_3v_4v_5$, without loss of generality, suppose that $v_1 v_5$ is red, then $v_2 v_4$ is blue, otherwise there is a red $2K_2$. Similarly, we have $v_1 v_3$ and $v_3 v_5$ are both red, thus $v_1 v_4$ and $v_2 v_5$ are both blue, yielding a blue $2K_2$, a contradiction.

If $n\geq m\geq 2$ and $n \geq 3$, there must be a vertex $u$ in $G$ with a red neighbor $u_{1}$ and a blue neighbor $u_{2}$. Then $V(G) \setminus\{u, u_{1}, u_{2} \}$ must induce a red $(m-1)K_2$ or a blue $(n-1)K_2$ by induction. Along with edge $u u_{1}$ or $u u_{2}$, we get a red $mK_2$ or a blue $nK_2$.
\hfill $\square$

\end{document}